\documentclass[letterpaper, 10 pt, conference]{ieeeconf}  
\usepackage{ifpdf}
\ifpdf
\usepackage[pdftex]{graphicx}
\usepackage{epstopdf}
\usepackage[pdftex,bookmarks=false,colorlinks=true,pdfstartview=FitBV,linkcolor=blue,citecolor=blue,urlcolor=blue]{hyperref}
\else
\usepackage{graphicx}
\usepackage{epsfig,psfrag}
\psdraft\psfull
\fi

\IEEEoverridecommandlockouts                              
\usepackage{amsmath,amssymb,amsfonts}
\usepackage{dsfont}
\usepackage[noend]{algpseudocode}
\usepackage{algorithm}
\usepackage{graphicx}
\usepackage{float}
\usepackage{subfigure}
\usepackage{setspace}
\usepackage{color}
\usepackage{cite}
\usepackage{romannum}
\usepackage{tabularx}
\usepackage{comment}
\usepackage{mathrsfs}
\usepackage{url}
\usepackage{enumerate}

\newcolumntype{C}[1]{>{\centering\arraybackslash}p{#1}}

\newcommand{\xe}{x_{\mbox{\tiny $E$}}}
\newcommand{\ye}{y_{\mbox{\tiny $E$}}}
\newcommand{\Xe}{p_{\mbox{\tiny $E$}}}

\newcommand{\thetae}{\theta_{\mbox{\tiny $E$}}}

\newcommand{\Eset}{\mathcal{E}}
\newcommand{\Pset}{\mathcal{P}}

\renewcommand{\cite}[1]{[\citen{#1}]}

\newtheorem{theorem}{Theorem}[section]
\newtheorem{lemma}[theorem]{Lemma}

\newtheorem{remark}{Remark}

\newtheorem{definition}[theorem]{Definition}
\newtheorem{conjec}[theorem]{Conjecture}

\renewcommand{\baselinestretch}{0.98}

\title{\LARGE \bf
  Apollonius Allocation Algorithm for Heterogeneous Pursuers\\ to Capture Multiple Evaders}

\author{Venkata Ramana Makkapati and Panagiotis Tsiotras%
    \thanks{This work has been supported by NSF award CMMI-1662542.}%
  \thanks{The authors are with the School of Aerospace Engineering, Georgia Institute of Technology, Atlanta. GA 30332-0150. USA.
    Email: {\small \{mvramana, tsiotras\}@gatech.edu}}%
}

\begin{document}

\raggedbottom

\maketitle
\thispagestyle{empty}
\pagestyle{empty}

\begin{abstract}
In this paper, we address pursuit-evasion problems involving multiple pursuers and multiple evaders.
The pursuer and the evader teams are assumed to be heterogeneous, in the sense that each team has agents with different speed capabilities.
The pursuers are all assumed to be following a constant bearing strategy.
A \emph{dynamic divide and conquer} approach, where at every time instant each evader is assigned to a set of pursuers based on the instantaneous positions of all the players, is introduced to solve the multi-agent pursuit problem.
In this regard, the corresponding multi-pursuer single-evader problem is analyzed first.
Assuming that the evader can follow any strategy, a dynamic task allocation algorithm is proposed for the pursuers.
The algorithm is based on the well-known Apollonius circle and allows the pursuers to allocate their resources in an intelligent manner while guaranteeing the capture of the evader in minimum time.
The proposed algorithm is then extended to assign pursuers in multi-evader settings that is proven to capture all the evaders in finite time.
\end{abstract}


\section{Introduction} \label{sec:intro}

Coordination strategies for unmanned aerial vehicles (UAVs)
has been an active area of research especially in the realm of multi-agent systems over the past decade \cite{Huang2013, Nex2014, Yuan2015}, having
numerous applications, including agriculture, aerial surveying, fire detection, disaster management, weather monitoring, and commercial product delivery.
Recent analyses of the commercial UAV market show that their use is expected to grow manyfold over the coming years,
as aerial drones are becoming a household product, used for recreational and industrial purposes alike\footnote{\url{http://www.businessinsider.com/commercial-uav-market-analysis-2017-8}}.
These advancements suggest an urgent need to explore designs for airspace safety systems that can regulate the traffic and usage of UAVs in a large scale.
Similarly, UAVs already play a major role in military engagement scenarios, and their use as part of swarm tactics (encirclement, coordinated attack, search and rescue,
perimeter defense) promises to change future battlefield operations.
It should therefore come as no surprise that a great amount of work has been devoted over the past decade to study coordination strategies of multi-agent UAV problems \cite{Rizk2018}.
To this end, UAV coordination strategies that formulate the problem as
a multi-player pursuit-evasion (PE) game offer solutions that address many of the challenges involving multi-agent systems such as of collision avoidance, surveillance and target acquisition~\cite{Mylvaganam2014, LasFargeas2015}.

The literature for multi-pursuer multi-evader (MPME) problems is actually limited.
In most cases, some form of heuristic is introduced in order to make the problem tractable.
Ge et al.~\cite{Ge2006} proposed three approaches, which include hierarchical decomposition, moving horizon hierarchical decomposition, and cooperative control.
Li et al.~\cite{Li2005} also explored a hierarchical approach, while Jin and Qu \cite{Jin2011} proposed a heuristic task allocation algorithm.
Extensions to the MPME problem includes problems with incomplete information \cite{Antoniades2003}, nonlinear dynamics \cite{Stipanovic2010}, and a mix of continuous and discrete observations \cite{Stipanovic2009}.
However, finding scalable algorithms which can be implemented in real-life MPME scenarios is still an open problem~\cite{Vieira2009a, Pierson2017}.

This work aims at extending current solution techniques for MPME games involving large teams of UAVs by
developing implementable, scalable solutions
based on a decomposition of the original MPME problem to a sequence of simpler multi-pursuer single-evader (MPSE) problems.
A major enabler for this decomposition is a new result that allows us to characterize each pursuer as \emph{active} or \emph{redundant} for each evader.
Only the relevant pursuers participate in the MPSE pursuit of each evader.
The identification and classification of each pursuer as active or redundant makes use of the classical tool of the Apollonius circle \cite{isaacs1999differential}.
Previously, the proposed MPME formulation was analyzed with the pursuit and the evading teams being homogeneous \cite{makkapati2019optimal}.
In this paper, we generalize the results in Ref. \cite{makkapati2019optimal} to include heterogeneous teams of agents with guarantees on finite-time capture.
The videos for the simulation results discussed in this paper can be found on the web\footnote{\url{https://youtu.be/KTH9lmdUdRs}}.

\section{Motivation and Problem Setup} 

Consider a group of $n$ agents (pursuers) guarding a given area of interest.
The objective of the agents is to pursue and intercept  $m$ (where typically $m \le n$) intruders (or evaders) that may be detected in this area.
Some of the relevant questions that arise while solving this problem include: a)Which pursuer(s) should go after which evader(s)?
b) How many pursuers should chase each intruder (evader) in order to capture it in the shortest time possible?
  
Obtaining the answers to the previous questions in their most general form is elusive at this point.
Addressing them involves solving a multi-player dynamic game, eventually demanding the solution to a high-dimensional partial differential equation, with the dimensionality increasing as the number of players ($n+m$).
In order to proceed and mitigate this problem, the following assumptions are made in this work.
\begin{enumerate}
  \item[A1:] The pursuers are all faster compared to the evaders.
  
  \item[A2:] The pursuers know the instantaneous positions and the velocities of all the evaders, and each pursuer follows a constant bearing strategy with respect to the evader to which it is assigned.
\end{enumerate}
The rationale behind these assumptions is as follows.
Under assumption A1, a constant bearing strategy guarantees capture.
Furthermore, constant bearing strategy has been implemented successfully in various aerial defense systems \cite{Shneydor}.

Formally, we address a pursuit-evasion problem in the Euclidean plane that involves $n$ pursuers and $m$ evaders.
The pursuers' objective is to capture all the evaders.
Capture occurs when one or more pursuers enter the capture zone of an evader (assumed here to be a disk of radius $\epsilon > 0$ centered at the instantaneous position of the evader).
Let
$P_i$ denote the $i^{th}$ pursuer and let $E_j$ denote the $j^{th}$ evader.
Let also $\Pset = \{P_1,P_2,\ldots,P_n\}$ denote the set of pursuers and, similarly, let $\Eset = \{E_1,E_2,\ldots,E_m\}$ denote the set of evaders.
The equations of motion of all the agents are given below
\begin{align}
\dot{x}_i &= u_i\cos \theta_i, \quad \dot{y}_i = u_i\sin \theta_i, \quad i \in  \Pset, \label{eqn:p1} \\
\dot{x}_j &= v_j\cos \theta_j, \quad \dot{y}_j = v_j\sin \theta_j, \quad j \in \Eset, \label{eqn:e}
\end{align}
where $p_i = (x_i,y_i) \in \mathbb{R}^2$, and $e_j = (x_j, y_j) \in \mathbb{R}^2$ denote the positions of pursuer $P_i$, and
evader $E_j$, respectively, and $\theta_i$ and $\theta_j$ denote the heading angles (control inputs) for the pursuers and the evaders, respectively.
In (\ref{eqn:p1}) and (\ref{eqn:e}), $u_i$ and $v_j$ are the speeds of $P_i$ and $E_j$, which are assumed to be constant with $\underset{i \in \mathcal{P}}{\min}~u_i>\underset{j \in \mathcal{E}}{\max}~v_j$ (A1).

A potential approach  to solve complicated MPME problems is a \emph{dynamic} ``divide and conquer'' approach,
where the pursuers are divided into several groups corresponding to the evader they pursue at each instant of time.
In essence, such  divide and conquer strategies formulate the original MPME problem as
a sequence of several (simpler) MPSE problems~\cite{Ge2006}.
This approach leads to decentralized (although likely suboptimal) solutions.
By analyzing the associated MPSE problems, one may arrive at an efficient dynamic task-allocation algorithm of pursuers to evaders.
To this end, we first present an approach for task allocation in MPSE problems using Apollonius circles in Section \ref{sec:strategies}.

\begin{figure}[htb!]
  \centering
  \includegraphics[width=0.4\textwidth]{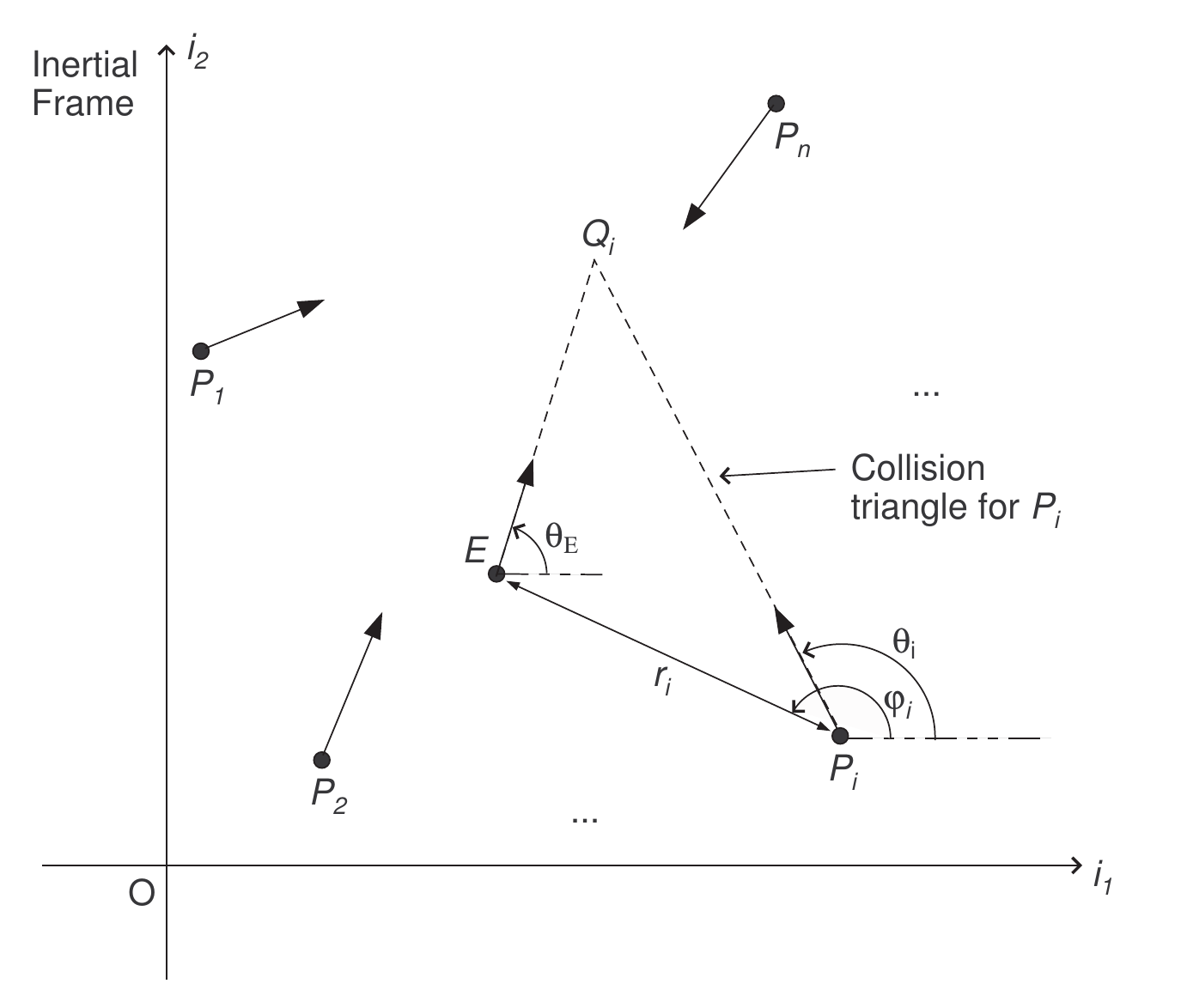}
  \caption{Schematic of the proposed multi-pursuer single-evader problem.}
  \label{fig:game_scene}
\end{figure}

A schematic of the proposed PE problem with one evader and multiple pursuers, that follow a constant bearing strategy--henceforth referred to as the MPSE problem--is shown in
Fig.~\ref{fig:game_scene}.
Since the pursuers are assumed to be following a constant bearing strategy, the problem can be analyzed by tracking the relative distances between the pursuers and the evader.
In this regard, the dynamics can be written in the form,
\begin{align}
\dot{r}_i &= v \cos(\thetae - \varphi_i) - u_i \cos(\theta_i - \varphi_{i}), \qquad i \in \Pset,
\end{align}
where $r_i$ is the relative distance between pursuer $P_i$ and the evader, and $\varphi_{i} = \text{atan2}(\ye-y_i,\xe-x_i)$ is the corresponding line of sight (LoS) angle.
For the MPSE problem, we drop the subscripts for the evader's speed $v$, and will use $E$ instead of $j$ (when required) to denote the single evader in this setting.
Furthermore, in the MPSE problem, we indicate the pursuers using the subscripts directly and the set $\mathcal{P} = \{1,2,\dots,n\}$.
Note that in the case of a constant bearing strategy the bearing angle between a pursuer and the evader remains constant until the time of capture.
Using this fact, the instantaneous heading of pursuer $P_i$ ($\theta_i$, $i \in \Pset$) can be obtained from the relation,
\begin{align}
u_i \sin(\theta_i - \varphi_{i}) &= v \sin(\thetae - \varphi_{i}), \label{eq:constri}
\end{align}
which is a function of the instantaneous heading of the evader $\thetae$. The above relation has two possible solutions for each $\theta_i$, given $\thetae$, and the solution for which $\dot{r}_i < 0$ is chosen.


\section{Task Allocation in Multi-Pursuer Single-Evader Problems}\label{sec:strategies}

In this section, we address task allocation in the MPSE problem by dynamically categorizing the pursuers into active and redundant.
The proofs for the lemmas presented in this section can be found in Ref. \cite{makkapati2019optimal}.
First, the formal definitions for time-of-capture and the corresponding capturing pursuer set, that are used to define active and redundant pursuers, in the MPSE problem are provided below.

\vspace*{1ex}

\begin{definition} \label{def:tc}
  For a given evading strategy, \emph{the time-to-capture $t_c~(\geq 0)$} is the minimum time so that there is at least one pursuer in the capture zone of the evader, and \emph{the capturing pursuer set $\mathcal{P}_c \subset \Pset$} is the set of pursuers that are in the capture zone of the evader at $t_c$.
\end{definition}

\vspace*{1ex}

The time-to-capture $t_c$ is always finite since the pursuers follow a constant bearing strategy.
Note that at the time of capture, one or more pursuers can be in the capture zone of the evader.
Therefore, $1 \leq \text{card}[\mathcal{P}] \leq n$, where card$[\cdot]$ represents the cardinality of the set.
Now, the following definitions establish the notions of active and redundant pursuers.

\vspace*{1ex}

\begin{definition}
  If there exists an evading strategy for which $P_i \in \mathcal{P}_c$, then $P_i$ is an \emph{active pursuer}.
  Otherwise, $P_i$ is a \emph{redundant pursuer}.
\end{definition}

\vspace*{1ex}

\begin{figure}[htb!]
  \centering
  \includegraphics[width=0.4\textwidth]{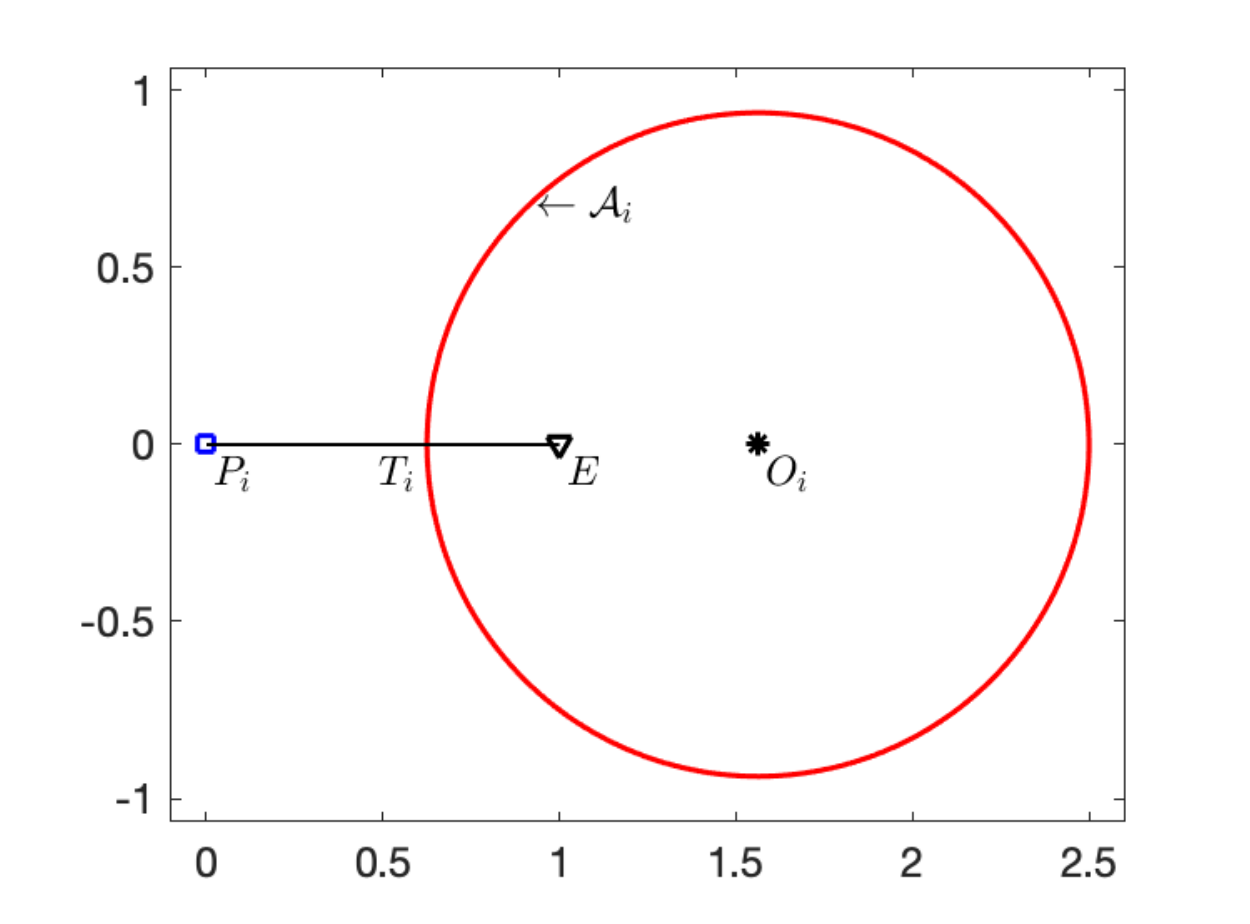}
  \caption{Apollonius circle. Simulation parameters: $u_i = 1$, $v=0.6$, $p_i(0) = (0,0)$, $\Xe = (1,0)$.}
  \label{fig:curves}
\end{figure}

Given the instantaneous positions of the pursuers and the evader, it is of interest to find a condition to verify whether a pursuer is active or redundant.
To this end, Apollonius circles are employed in this work.
The Apollonius circle for a pursuer-evader pair is the locus of points where capture occurs, for all possible initial headings of a non-maneuvering evader, given the initial positions of the
pursuer/evader pair and assuming that the pursuer follows a constant bearing strategy, see Fig.~\ref{fig:curves}.
For the MPSE problem, the Apollonius circle of the pair $P_i$-$E_j$ is denoted as $\mathcal{A}_{ij}$.
It has  its center at $O_i\bigg(\dfrac{x_j - \rho_{ij} x_i}{1 - \rho_{ij}^2},\dfrac{y_j - \rho_{ij} y_i}{1 - \rho_{ij}^2} \bigg)$ and radius $d_{ij}=\dfrac{\rho_{ij}\|p_i-p_j\|}{1-\rho_{ij}^2}$, where $\rho_{ij} = v_j/u_i$ (speed ratio) \cite{Ramana2017a}.
The Apollonius circles evolve in time as the players move, but the time dependencies will be dropped for the sake of brevity.
Let $T_{ij}$ be the closest point to the evader on the Apollonius circle where collision occurs when the evader goes head-on with the pursuer, as shown in Fig.~\ref{fig:curves}.
Therefore the distance of $T_{ij}$ from the evader is ${v_j\|p_i - p_j\|}/{(u_i+v_j)}$.
Next, we define \emph{Apollonius boundary} and analyze its properties to identify the active and the redundant pursuers.

\begin{definition} \label{def:AB}
  Given the positions all the pursuers and the evader $E_j$ at time $t \geq 0$,
  the \emph{Apollonius boundary} around the evader $E_j$ at time $t$ is the set of points 
  \begin{align*}
      \mathcal{B}^t_j = \{ X \in {\textstyle \cup_{i=1}^{n}} \mathcal{A}_{ij} ~|~ \mathcal{M}(e_j, X) \cap \big( {\textstyle \cup_{i=1}^{n}} \mathcal{A}_{ij} \big) = \{X\} \},
  \end{align*}
  where $\mathcal{M}(e_j,X)$ denotes the set of points on the line segment with endpoints $e_j$ (position of $E_j$) and $X$ at time $t$.
\end{definition}

\vspace*{1ex}

In other words, the Apollonius boundary is the set of points that belong to the union of all the instantaneous Apollonius circles corresponding to the evaders and, in addition,
each such point is the closest to the evader along its respective line-of-sight originating from the evader.
Note that the Apollonius boundary evolves with time as the Apollonius circles also evolve with time.
The subscript $j$ will be dropped for the rest of this section, since we are dealing with the MPSE problem.
\begin{figure}[htb!]
  \centering
  \includegraphics[width=0.4\textwidth]{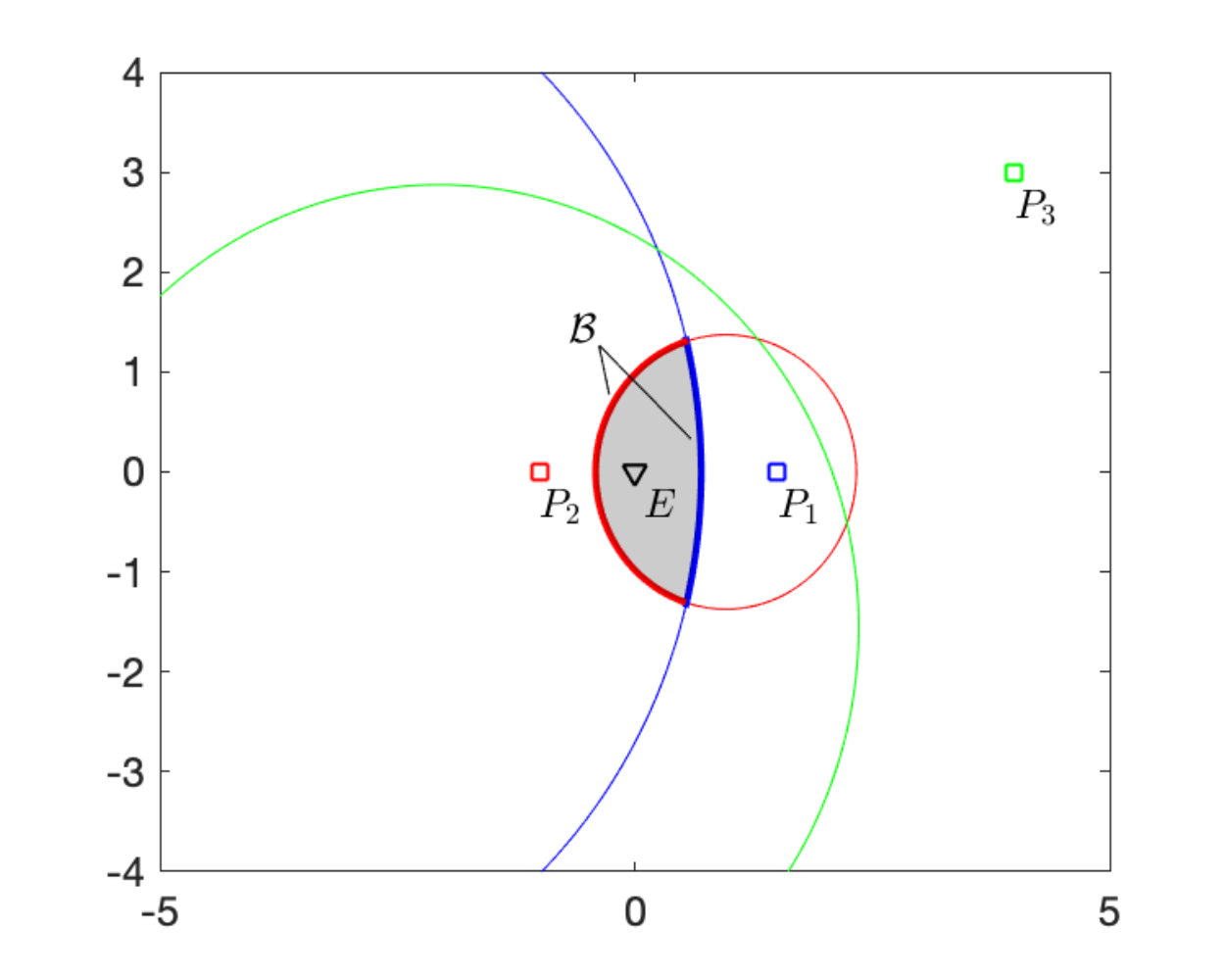}
  \caption{Apollonius boundary. Simulation parameters: $u_1 = 0.8$, $u_2 = 1$, $u_3 = 1.2$, $v=0.7$}
  \label{fig:boundaries}
\end{figure}

\subsection{Apollonius circle based Active Pursuer Check}

The algorithm that is developed to identify active/redundant pursuers in this work is based on the following conjecture.

\vspace*{1ex}

\begin{conjec} \label{conjecture:active}
  Given the positions of all the players in an MPSE problem at time $0 \leq t < t_c$, and assuming that the pursuers follow a constant bearing strategy, pursuer $P_i$ is active at time $t$ if $\mathcal{B}^t \cap \mathcal{A}_i \neq \varnothing$, and is redundant otherwise.
\end{conjec}

\vspace*{1ex}

The conjecture implies that a pursuer is active at time $0 \leq t < t_c$ if and only if its corresponding Apollonius circle is part of the Apollonius boundary at that instant.
The conjecture is inspired from the fact that the region in which the capture point lies in is bounded by the instantaneous Apollonius circle for any strategy of the evader.
Note that if a pursuer is active at time $t'$, it need not remain active for all $t > t'$.
But if a pursuer is redundant at time $t'$, it will remain redundant for all $t > t'$.
The following lemmas based on Conjecture~\ref{conjecture:active} provide simple checks to determine whether a pursuer is active or redundant.

\vspace*{1ex}

\begin{lemma} \cite{makkapati2019optimal} \label{lemma:one}
  Given the positions of the players in the MPSE problem at time $0 \leq t < t_c$, pursuer $P_i$ is the only active pursuer if and only if the conditions
  \begin{align}
  \mathcal{A}_i \cap \big( {\textstyle \cup_{j=1,\,j \neq i}^{n}} \mathcal{A}_j \big) &= \varnothing, \label{eq:lemma1_cond1}\\
  \mathcal{M}(e,T_i) \cap \big( {\textstyle \cup_{j=1,\,j \neq i}^{n}} \mathcal{A}_j \big) &= \varnothing, \label{eq:lemma1_cond2}
  \end{align}
  are satisfied, where $T_i$ is the closest point to the evader on the Apollonius circle $\mathcal{A}_i$.
  Furthermore, if conditions (\ref{eq:lemma1_cond1}) and (\ref{eq:lemma1_cond2}) are not satisfied, then $P_i$ is a redundant pursuer.
\end{lemma}

\vspace*{1ex}

\begin{lemma} \cite{makkapati2019optimal} \label{lemma:two}
  Given the positions of the players in the MPSE problem at time $0 \leq t < t_c$, and that the Apollonius circle $\mathcal{A}_i$ intersects at least one  of the other Apollonius circles.
  Then, pursuer $P_i$ is an active pursuer if and only if there exists $X \in \mathcal{I}_i$ such that $\mathcal{M}(e,X) \cap \big( {\textstyle \cup_{j=1}^{n}} \mathcal{A}_j \big) = \{X\}$,
  where $\mathcal{I}_i$ is the set of intersection points between $\mathcal{A}_i$ and the rest of the Apollonius circles.
\end{lemma}

\vspace*{1ex}

The set of intersection points $\mathcal{I}_i$ can be obtained analytically given the instantaneous positions of all the players \cite{Weisstein}.
The above two lemmas can be used to verify if a pursuer is active or redundant.
In this regard, Algorithm~\ref{algo:active} below, named Apollonius circle based Active Pursuer Check (AAPC), can be employed to check the status of each pursuer.
The time complexity of the algorithm is of order $O(n^2)$, since the maximum number of intersections between any two circles is  two.
Note that by dynamically allocating the task of capturing the evader using AAPC (where at every instant the active pursuers keep pursuing the evader while the redundant pursuers do not react), the pursuers as a group will be able to capture the evader in minimum time.
Furthermore, if a pursuer becomes redundant at any point of time $0 \leq t < t_c$, it remains redundant after that (i.e., till capture occurs).

\begin{algorithm}
  \caption{Apollonius circle based Active Pursuer Check (AAPC)}\label{algo:active}
  \begin{algorithmic}[1]
    \Require Positions of all the players ($p_1$,\dots,$p_n$,$e$,$i$)
    \Ensure Status of pursuer $P_i$
    \Procedure{obtain\_status}{$p_1$,\dots,$p_n$,$e$,$i$}
    \State intersection = 0 
    \State status = redundant
    \For  {$j=1$ to $n$ and $j \neq i$}
    \State Obtain the set of intersection points $\mathcal{I}_{ij}$
    \If {$\mathcal{I}_{ij} \neq \varnothing$}
    \State intersection = 1
    \For  {$\ell=1$ to $\text{card}[\mathcal{I}_{ij}]$}
    \State boundary = 1
    \For  {$k=1$ to $n$ and $k \neq i,\,j$}
    \If {$\mathcal{M}(e,X_\ell)$ intersects $\mathcal{A}_k$}
    \State boundary = 0
    \EndIf
    \EndFor
    \If {boundary = 1}
    \State status = active
    \State break from outermost loop.
    \EndIf
    \EndFor
    \EndIf
    \EndFor
    \If {intersection = 0}
    \State status = active
    \For {$j=1$ to $n$ and $j \neq i$}
    \If {$\mathcal{M}(e,T_i)$ intersects $\mathcal{A}_j$}
    \State status = redundant
    \State \textbf{break}
    \EndIf
    \EndFor
    \EndIf
    \State \Return status
    \EndProcedure
  \end{algorithmic}
\end{algorithm}

\subsection{Numerical Simulations} \label{subsec:numsim}

In this section simulations of pursuer allocation using AAPC involving five pursuers and one evader are presented.
There are three different pursuers with speeds $u_i \in \{0.8,1,1.2\}$, and the speed of the evader is set to
$v = 0.6$.
The radius of capture is chosen as $\epsilon = 0.1$.
The evader follows a form of blind evasion strategy with switching times that are predefined~\cite{Morgan2016}.
At each switching time the evader randomly chooses a heading from a set of allowable headings.
The allowable headings set that is specific to the example showcased in this work is $\{-\pi/4,~\pi/2,~3\pi/4\}$.

\begin{figure}[htb!]
  \centering
  \subfigure[Initial Apollonius circles]{\includegraphics[width=0.39\textwidth]{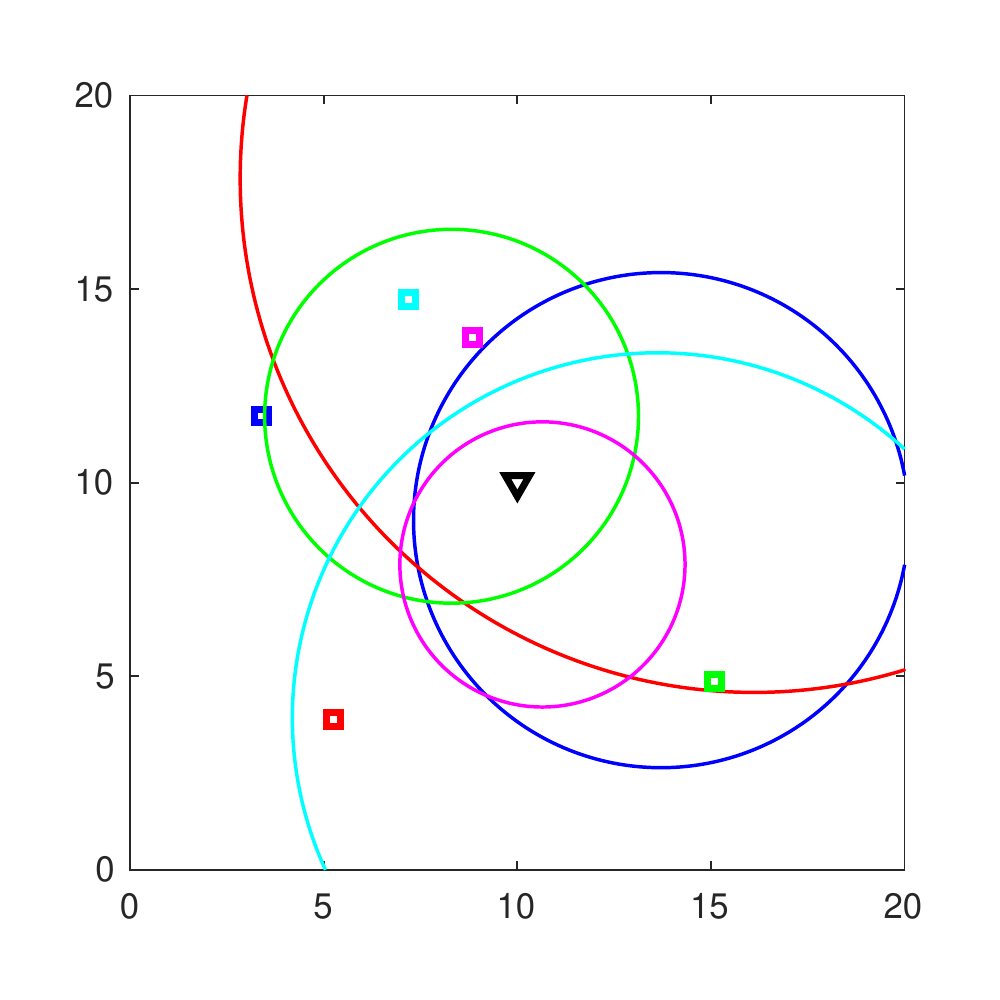}}  \quad
  \subfigure[Trajectories]{\includegraphics[width=0.39\textwidth]{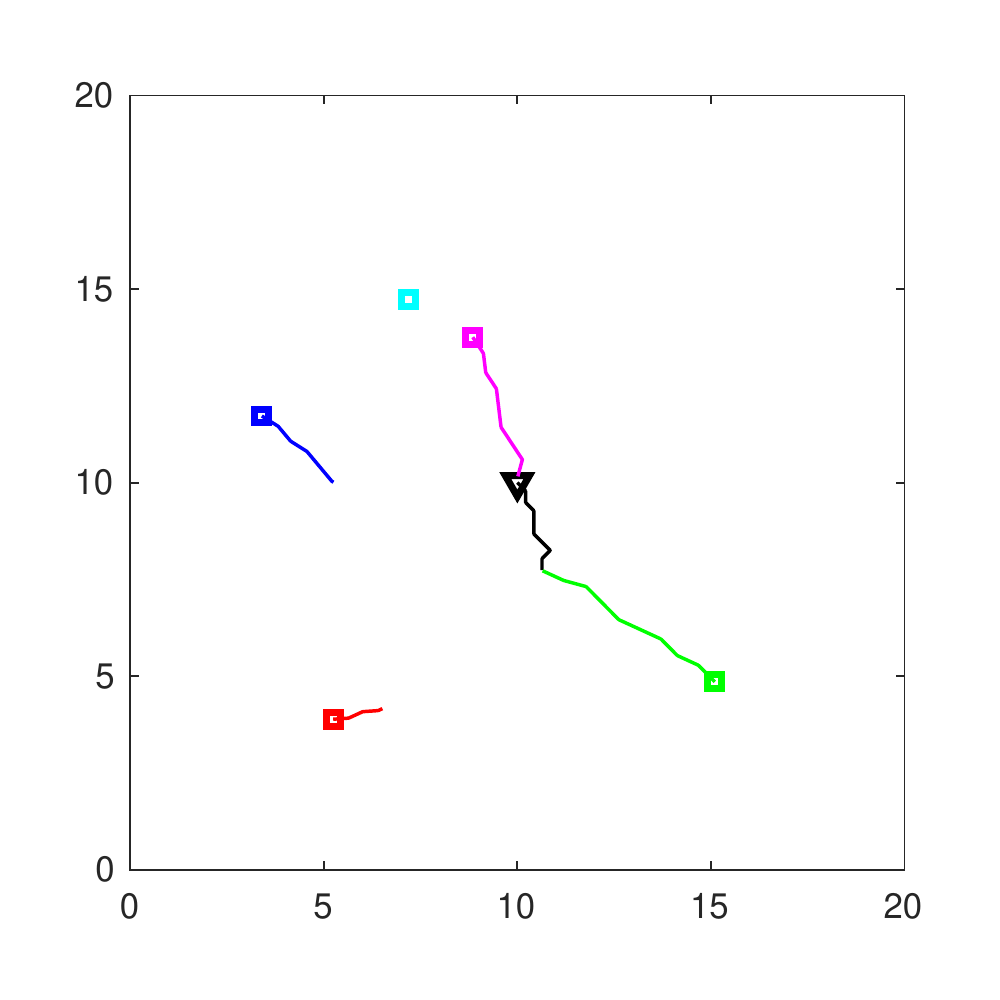}}
  \caption{Results obtained using AAPC for task allocation in the case of CB.}
  \label{fig:sim_cb}
\end{figure}

Fig.~\ref{fig:sim_cb}(a) shows the initial positions of all the players along with the corresponding Apollonius circles.
The triangle denotes the initial position of the evader and the square markers denote the initial positions of the pursuers.
It can be observed that at the initial time, the pursuers identified with the colors red, magenta, green, and blue are the active pursuers, as their corresponding Apollonius circles are part of the Apollonius boundary.
Fig.~\ref{fig:sim_cb}(b) shows the trajectories of all the players.
It can be seen that the green pursuer finally captures the evader, and the rest of the three pursuers, which are  initially active, become redundant as time progresses.
The cyan pursuer, which is redundant at the initial time, does not move at all.


\section{Apollonius Allocation Algorithm}\label{sec:mpme}

In this section, the AAPC is extended to solve MPME problems.
Given the positions of all the players at some instant of time, the set of evaders for which a pursuer is active can be obtained using AAPC.
Note that at a given time instant, a pursuer can be classified as active by more than one evader or no evader whatsoever.
In the case where a pursuer is classified as active for more than one evader, one can break the tie by assigning the pursuer to the evader that can be captured in minimum possible time i.e., $\mathrm{argmin}_j \|p_i - e_j\|/(u_i + v_j)$.

\vspace*{1ex}

\begin{remark}
    The aforementioned criterion for breaking a tie is equivalent to choosing the nearest evader when the teams are assumed to be homogeneous \cite{makkapati2019optimal}.
\end{remark}

\vspace*{1ex}

Using this criterion, the following algorithm can be used for pursuer allocation in MPME problems.

\vspace{1em}

\emph{Apollonius Allocation (A2) Algorithm:}
At a given time instant $0 \leq t \leq t_c$, let $\Eset_f$ be the set of evaders that are yet to be captured, and let $\Eset_c$ be the set of evaders that have already been  captured.
Note that $\mathcal{E} = \Eset_f \cup \Eset_c$.
Given the current positions of all the players,
let $\mathscr{I}: \Eset_f \rightarrow 2^{\Pset}$ be the initial allocation function that maps each evader $E_j$ (in $\Eset_f$) to its set of active pursuers obtained by considering the positions of all the pursuers.
That is, for a given $j \in \Eset_f$, $\mathscr{I}(j)$ is a subset of $\Pset$.
Furthermore, $\Pset_a = \cup_{j \in \Eset_f}\mathscr{I}(j)$ denotes the set of all the active (or assigned) pursuers according to the
initial allocation function $\mathscr{I}$.
Given the initial allocation function $\mathscr{I}$, let now $\mathscr{J}: \Pset \rightarrow 2^{\Eset_f}$ be the dual function defined by
$\mathscr{J}(i) = \{ j \in \Eset_f: \mathscr{I}(j) = i \}$.
In other words, $\mathscr{J}$ maps each pursuer to the set of the evaders to which it is allotted as per $\mathscr{I}$.
Next, we define the \emph{final allocation function} $\mathscr{F}$ and the \emph{intermediate allocation function} $\mathscr{G}$
as follows.

\begin{enumerate}[(a)]
  
  \item
  If card$[\mathscr{J}(i)] \leq 1$, for all $i  \in \Pset$, then let $\mathscr{F} = \mathscr{I}$.
  Otherwise, let $\mathscr{G}: \Eset \rightarrow 2^{\Pset}$ be defined as $\mathscr{G}(j) = \Big\{ i \in \mathscr{I}(j) : j = \underset{k \in \mathscr{J}(i)}{\mathrm{argmin}}~\dfrac{\|p_i - e_k\|}{u_i + v_k} \Big\}$.
  The function $\mathscr{G}$ maps each evader to a set of pursuers in accordance to the mapping $\mathscr{I}$,
  such that each active pursuer is assigned to the nearest evader among its assigned ones.
  Note that $\mathscr{G}(j)$ can be an empty set for some $j$, i.e., an evader can end up be unassigned as per $\mathscr{G}$.

  \item
  Let $\Pset_u  = \Pset \backslash \Pset_a$ be the set of unassigned pursuers.
  Now for each evader $E_j$, find the active pursuers considering the positions of the pursuers that are only in the set $\mathscr{G}(j)~\cup~\Pset_u$,
  and obtain an \emph{updated allocation function} $\mathscr{I}'$ and its corresponding dual $\mathscr{J}'$.
  
  \item Repeat steps (a) and (b), by replacing $\mathscr{I}$ and $\mathscr{J}$ with $\mathscr{I}'$ and $\mathscr{J'}$, respectively, until $\mathscr{F}$ is obtained.   
    
\end{enumerate}

Note that in step (a) of the algorithm, ties with multiple assignments of the same pursuer are broken using distance as the metric.
Furthermore, if each pursuer is assigned to only one evader or if it remains unassigned, then the initial allocation function $\mathscr{I}$ is also the final one.
In any other case, once the intermediate allocation function $\mathscr{G}$ is obtained in step (b) of the algorithm,
the set of unassigned pursuers according to $\mathscr{I}$ is obtained.
In step (b), an updated allocation function $\mathscr{I}'$ is obtained by checking for active pursuers among the set of unassigned pursuers coupled with the pursuers assigned as per $\mathscr{G}$, for each evader.
Because one of the unassigned pursuers (in the set $\Pset_u$) can become active to the evaders that have lost one or more pursuers during the tie break in step (a).
With $\mathscr{I}'$ and its corresponding dual $\mathscr{J}'$, steps (a) and (b) are repeated until each pursuer has only one (or none) assignment.
Once an evader is captured, it is removed from the set $\Eset_f$ and added to the set $\Eset_c$.

\begin{figure*}[htb!]
  \centering
  \subfigure[$t$ = 0 (At initial time)]{\includegraphics[width=0.32\textwidth]{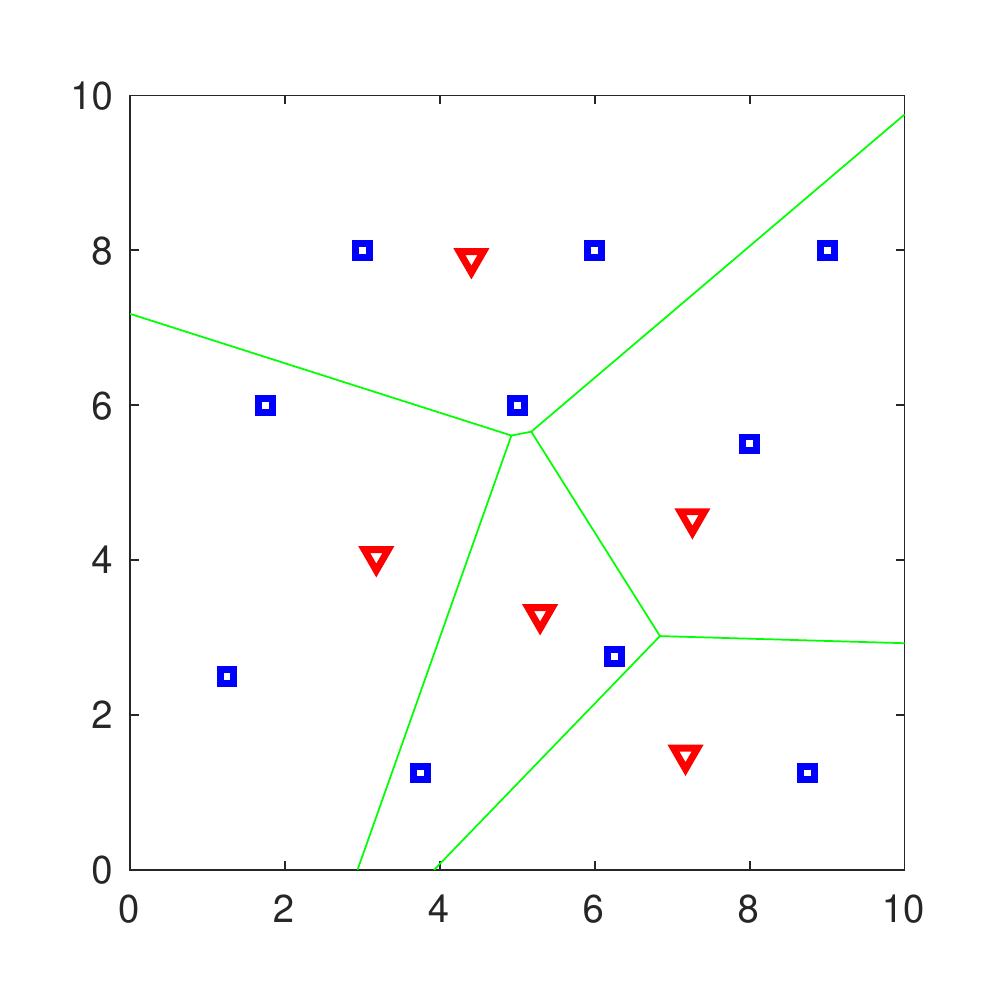}}
  \subfigure[$t$ = 1.8 (Three evaders are captured)]{\includegraphics[width=0.32\textwidth]{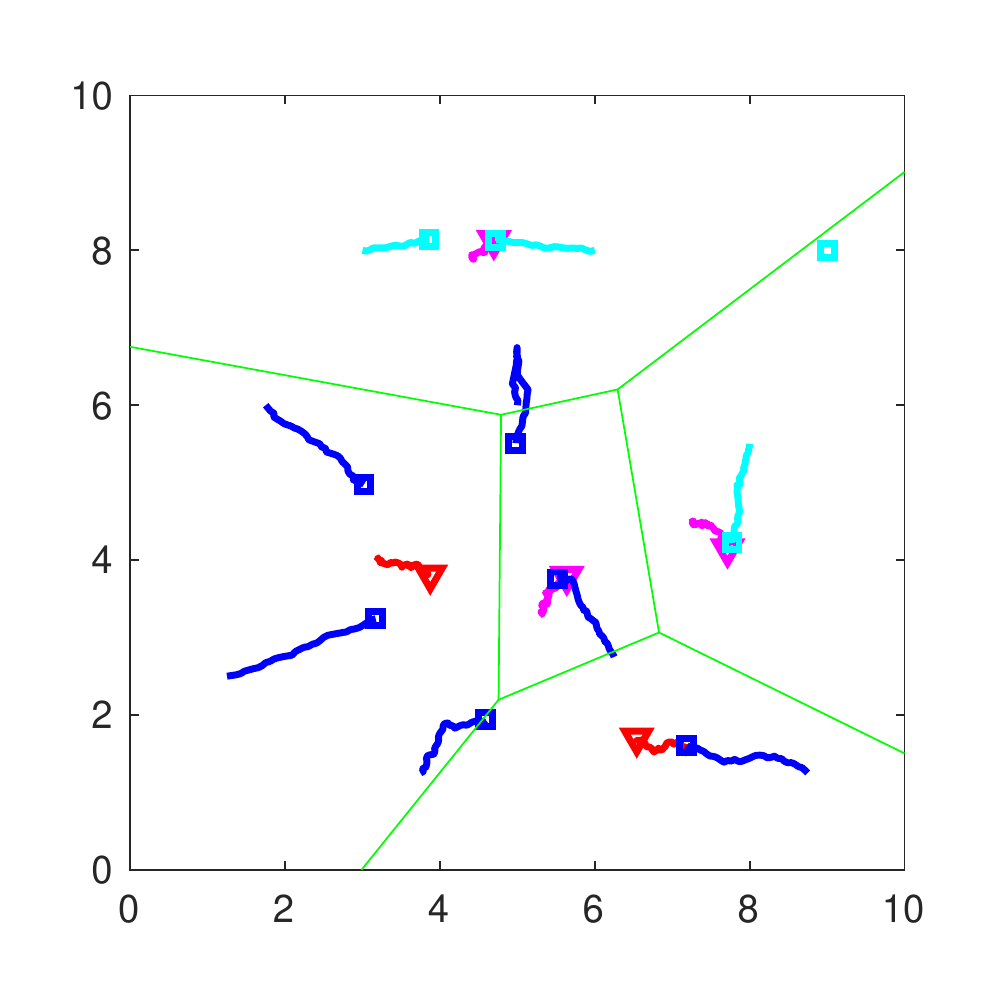}}
  \subfigure[$t$ = 2.75 (All evaders are captured)]{\includegraphics[width=0.32\textwidth]{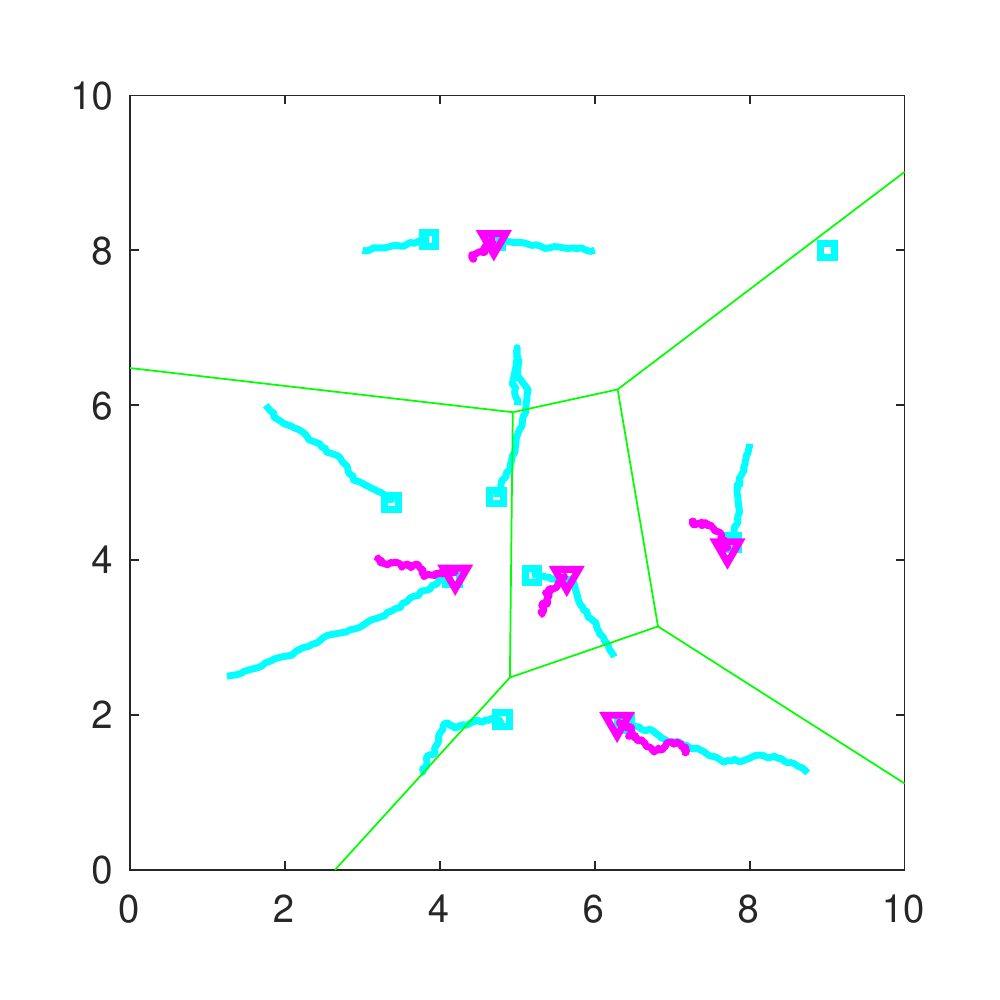}}
  \caption{Plots showing the positions and trajectories of the players in a multi-pursuer (squares) multi-evader (triangles) problem at different time instants.}
  \label{fig:sim_mpme}
\end{figure*}

The above algorithm is run at every time instant to obtain the allocation function $\mathscr{F}$, given the players' current positions, until all the evaders are captured, 
i.e., until $\mathcal{E}_f$ is empty.
The algorithm provides a potentially sub-optimal solution, but it is scalable for any number of pursuers and evaders.
The algorithm guarantees capture of all $m$ evaders as is shown in Theorem~\ref{Theorem:main} below.
In order to prove this theorem, several preparatory results are needed.

\vspace*{1ex}

\begin{definition}\label{def:CST}
  Given the instantaneous positions of the players in the MPME problem at time $t \geq 0$, the \emph{current shortest time} (CST) is defined by $t_s = \underset{(i,j) \in \Pset \times \Eset_f}{\min}~\dfrac{\|p_i - e_j\|}{u_i + v_j}$.
\end{definition}

\vspace*{1ex}

\begin{lemma}
  At a given time instant $t \geq 0$, $i^* \subseteq \mathscr{F}(j^*)$, where $(i^*,j^*) = \underset{(i,j) \in \Pset \times \Eset_f}{\mathrm{argmin}}~\dfrac{\|p_i - e_j\|}{u_i + v_j}$, and $\mathscr{F}$ is the final allocation function of A2 algorithm. \label{lemma:MPME_2}
\end{lemma}

\vspace*{1ex}

\begin{proof}
    From Definition~\ref{def:CST} it is understood that $t_s$ represents the minimum possible time taken to capture an agent in the evading team by an agent in the pursuing team, when the corresponding evader goes head-on with the pursuer.
    Therefore, if it can be shown that the Apollonius boundary around the evader $j^*$ contains part of the Apollonius circle $\mathcal{A}_{i^*j^*}$ (i.e., $\mathcal{B}^t_{j^*} \cap \mathcal{A}_{i^*j^*} \neq \varnothing$, $\forall
    ~t \geq 0$), then the pursuer $i^*$ will be assigned to the evader ${j^*}$, even when there is a tie.
    Since the Apollonius circle denotes the capture points for a non-maneuvering evader, the point $T_{i^*j^*}$, with the length of the line segment $E_{j^*}T_{i^*j^*} = v_{j^*}t_s$, is the closest capture point to evader $j^*$ along the corresponding line of sight.
    Therefore $T_{i^*j^*} \in \mathcal{B}^t_{j^*}$, and hence, $\mathcal{B}^t_{j^*} \cap \mathcal{A}_{i^*j^*} \neq \varnothing$, $\forall~t\geq 0$.
\end{proof}

\vspace*{1ex}

\begin{lemma}  \label{lemma:MPME_3}
  Assuming the pursuers are assigned to the evaders using A2, at any given time $t \geq 0$, CST will converge to zero in finite time, and hence at least one evader will be captured in finite time. 
\end{lemma}

\vspace*{1ex}

\begin{proof}
  From Lemma~\ref{lemma:MPME_2}, pursuer $i^*$ (corresponding to the CST) is always assigned to evader $j^*$. Since all the pursuers are faster compared to the evaders, and since they follow a constant bearing strategy, d$t_s/$dt $\leq (v_{j^*} - u_{i^*})/(v_{j^*} + u_{i^*}) < 0$, for all $t \geq 0$~\cite{Shneydor}. 
  Furthermore, as the initial CST is finite, the CST converges to zero in finite time. 
  Hence, capture of one evader is guaranteed in finite time.
\end{proof}

\vspace*{1ex}

\begin{theorem} \label{Theorem:main}
  The A2 algorithm guarantees capture of all the evaders in finite time.
\end{theorem}

\vspace*{1ex}

\begin{proof}
  The result immediately follows from Lemmas \ref{lemma:MPME_2} and \ref{lemma:MPME_3}. Note that the CST is updated (from zero) every time a capture occurs, and the captured evader is removed from the list of participating players. Also, the number of evaders are finite.
\end{proof}

\vspace*{1ex}

Fig.~\ref{fig:sim_mpme} demonstrates the performance of A2 algorithm for 10 pursuers and 5 evaders.
For the sake of simulation, the evading is assumed to be homogeneous.
The simulation parameters remain the same as in Section~\ref{subsec:numsim}.
In Fig.~\ref{fig:sim_mpme}, the red triangles indicate the current positions of the evaders that are not captured and the magenta ones are the evaders that are captured.
The blue squares indicate the current positions of the active pursuers and the cyan ones indicate the redundant pursuers.
In all three plots, the Voronoi partition of the domain with the evaders as generators is also included for reference.


\section{Conclusion} \label{sec:conclude}

Under the assumption that the pursuers are faster than the evader(s), and that they follow a constant bearing strategy, workable solutions for multi-pursuer single-evader (MPSE) and multi-pursuer multi-evader (MPME) problems are provided.
A dynamic allocation algorithm for the pursuers that is independent of the evader's strategy has been proposed to solve the MPSE problem.
The proposed algorithm is based on the notion of active/redundant pursuers, and employs the concept of Apollonius circles.
The algorithm is further extended to solve MPME problems for any number of pursuers and evaders.
These algorithms ensure capture of all the evaders either in an MPSE or an MPME setting in finite time.
Several extensions of this work are possible.
For example,
the computational requirements can be reduced by having an estimate of when the assignment can change to avoid unnecessary calculations at every time instant.
It would also be very interesting to extend
the notion of Apollonius circles to account for turn-radius constraints for all the players, but this would probably end up being a major challenge.


\vspace*{1ex}

\end{document}